\documentclass[12pt]{article}
\usepackage{latexsym}
\usepackage{amssymb}
\usepackage{amsmath}
\usepackage{rotating}
\usepackage{multirow}
\usepackage{mathrsfs}
\usepackage{wasysym}
\usepackage{esint}
\usepackage{rawfonts}
\usepackage{hyperref}

\input{prepictex}
\input{pictex}
\input{postpictex}
\usepackage[OT2,OT1]{fontenc}
\def\cyr{%
\renewcommand\rmdefault{wncyr}%
\renewcommand\sfdefault{wncyss}%
\renewcommand\encodingdefault{OT2}%
\normalfont
\selectfont}
\DeclareMathAlphabet{\zap}{OT1}{pzc}{m}{it}
\DeclareTextFontCommand{\textcyr}{\cyr}
\def\be{\begin{equation}}
\def\ee{\end{equation}}
\def\bea{\begin{eqnarray*}}
\def\eea{\end{eqnarray*}}

\newtheorem{main}{Theorem}
\DeclareMathOperator{\Hess}{Hess}

\newtheorem{thm}{Theorem}
\newtheorem{lem}{Lemma}
\newtheorem{prop}{Proposition}

\newenvironment{proof}{\medskip \noindent
{\bf Proof. }}{\hfill \rule{.5em}{1em}
\\}

\def\CP{{\mathbb C \mathbb P}}

\begin{document}
\author{Caner Koca}
\title{Extremal K\"ahler Metrics and Bach--Merkulov Equations}\date{}
\maketitle

\begin{abstract}
In this paper, we study a coupled system of equations on oriented compact 4-manifolds which we call the Bach--Merkulov equations. These equations can be thought of as the conformally invariant version of the classical Einstein--Maxwell equations. Inspired by the work of C. LeBrun on Einstein--Maxwell equations on compact K\"ahler surfaces, we give a variational characterization of solutions to Bach--Merkulov equations as critical points of the Weyl functional. We also show that extremal K\"ahler metrics are solutions to these equations, although, contrary to the Einstein--Maxwell analogue, they are not necessarily minimizers of the Weyl functional. We illustrate this phenomenon by studying the Calabi action on Hirzebruch surfaces.
\end{abstract}

\section{Introduction}
Let $M$ be a smooth oriented $n$-manifold. A Riemannian metric $g$ on $M$ is said to satisfy the \textit{Einstein--Maxwell equations} if
\begin{equation}\label{eq1}
\begin{array}{r}[r+F\circ F]_\circ = 0\\[10pt]
dF=0,\quad d*F = 0
\end{array}
\end{equation}
for some 2-form $F$ on M. Here, $r$ is the Ricci tensor of $g$; $(F\circ F)_{ij} = {F_i}^s F_{sj}$ is the composition of $F$ with itself as an endomorphism of the tangent bundle $TM$; $[\cdot]_\circ$ denotes the trace-free part of a $(2,0)$-tensor, and $*$ is the Hodge operator with respect to the metric $g$. When $M$ is \emph{compact}, the second line of \eqref{eq1}, which is called \emph{Maxwell equations}, is equivalent to saying that $F$ is \emph{harmonic} with respect to $g$, i.e. $\Delta F=0$.

By Hodge theory we know that any harmonic form $F$ minimizes the $L^2$ norm $F\mapsto\int_M |F|^2_g d\mu_g$ among the forms cohomologous to $F$, namely on $[F]\in H_{\mathrm{dR}}^2(M,\mathbb R)$. If, in addition, $M$ has dimension $4$, the integral $\int_M |F|_g d\mu_g$ is unchanged if $g$ is replaced by any conformally related metric $\tilde g := ug$, for a positive smooth function $u$ on $M$. Therefore, if $F$ is harmonic with respect to $g$, it will be harmonic with respect to $\tilde g$. By contrast, the first line of \eqref{eq1} is certainly not conformally invariant in any dimension. There is, however, an interesting conformally invariant counterpart of these equations introduced by Merkulov in \cite{M84}:
\begin{eqnarray}\label{eq2}
\begin{array}{rll}
B+ [F\circ F]_\circ = 0\\[10pt]
dF=0,\quad d*F = 0
\end{array}
\end{eqnarray}
where $B_{ij} = (\nabla^s \nabla^t +\frac{1}{2}r^{st}) W_{isjt}$  is the \textit{Bach tensor} \cite{B}. When $M$ is compact, this tensor arises as the Euler-Lagrange equations for the \emph{Weyl energy functional} $g\mapsto \int_M |W|^2 d\mu_g$ over the space of all metrics. That is, if we vary the metric $g_t=g_o+t h +o(t^2)$, then \cite{Bes}
\begin{equation}\label{bachvar}
\left.\frac{d}{dt}\right|_{t=0} \mathcal W (g_t)= \int_M \langle h,B\rangle d\mu = \int_M g^{is}g^{tj} h_{st}B_{ij} d\mu.
\end{equation}
Note that in 4 dimensions, $\mathcal W$ is indeed conformally invariant since the conformal change $\tilde{g} = u g$ of the metric implies
$$d\tilde\mu = u^2 d\mu \quad\textnormal{and}\quad {{{\tilde{W}}_{ijk}}}{}^l = {W_{ijk}}^l.$$
Bach tensor, too, behaves well under conformal change: ${\tilde B}_{ij}= \frac{1}{u}B_{ij}$. To see this note that for the rescaled variation $\tilde g_t = ug_0 + t uh + o(t^2)$ we have
$$\left.\frac{d}{dt}\right|_{t=0} \mathcal W (\tilde g_t)= \int_M \langle uh,\tilde B\rangle d\tilde\mu = \int_M \tilde g^{is}\tilde g^{tj} u h_{st} B_{ij} d\tilde\mu$$
and comparing it to \eqref{bachvar} we deduce that $\tilde B_{ij} = \frac{1}{u} B_{ij}$. Also $B$ is symmetric, trace-free and divergence-free. Note also that $[F\circ F]_\circ$, the other term in \eqref{eq2}, rescales similar to $B_{ij}$ under conformal rescaling. Clearing out the $\frac{1}{u}$ factors, we see that, when $M$ is a compact manifold of dimension $4$, the coupled system of equations \eqref{eq2} is conformally invariant in the sense that if $(g,F)$ is a solution, so is $(ug, F)$ for any positive smooth function $u$.

Both Einstein--Maxwell and Bach--Merkulov equations stem from a variational origin. For any given de Rham class $\Omega\in H^2_{\mathrm dR}(M,\mathbb R)$, solutions $(g,F)$ of Einstein-Maxwell equations with $[F]=\Omega$ are in fact the critical point of the coupled action
\begin{align*}
\mathcal G_1 \times \Omega & \longrightarrow \mathbb R \\
(g,F) &\longmapsto \int_{M} s_g + |F|_g^2 d\mu_g
\end{align*}
where $\mathcal G_1$ stands for the space of unit volume metrics \cite{L08}. Similarly \cite{B},\cite{M84}, Bach--Merkulov equations are the critical point of the action
\begin{align*}
\mathcal G_1 \times \Omega & \longrightarrow \mathbb R \\
(g,F) &\longmapsto \int_{M} |W|^2_g + |F|_g^2 d\mu_g.
\end{align*}

In \cite{L08}, C. LeBrun studied Einstein--Maxwell equations \eqref{eq1} on compact smooth 4-manifolds, and discovered some fascinating properties of these equations in relation to K\"ahler geometry. He showed that constant scalar curvature K\"ahler metrics satisfy \eqref{eq1}; all solutions to \eqref{eq1} are critical points of $L^2$-norm of scalar curvature on $\mathcal G_\Omega$, the space of metrics for which a fixed cohomology class $\Omega$ is represented by a self-dual harmonic form $\Omega_g$; and on complex surfaces constant scalar curvature K\"ahler metrics are global minimizers of that action if $c_1\cdot\Omega\leq 0$. These results are summarized in Section 2. The aim of this paper is to state and prove the relevant properties of the Bach--Merkulov equations. The main results are:

\begin{main} \label{thmA}
Let $M$ be a compact complex surface, and let $g$ be a metric conformal to an extremal K\"ahler metric on $M$. Then $g$ solves the Bach--Merkulov equations for some $F$. As a consequence, on any compact complex surface K\"ahler type we can solve \eqref{eq2}.
\end{main}

In other words, extremal K\"ahler metrics are standard solutions on a compact complex surface. On a more general $4$--manifold the Bach--Merkulov equations naturally become critical points of the Weyl energy functional $g\mapsto\int_M |W^+|^2 dvol$:

\begin{main}\label{thmB}
Let $M$ be a smooth compact oriented $4$--manifold, and $\Omega\in H_{\mathrm dR}^2(M,\mathbb R)$ be any de Rham class. A metric $g\in\mathcal G_\Omega$ is a critical point of the restriction of Weyl functional to $\mathcal G_\Omega$ iff $g$ is a solution of Bach--Merkulov equations in conjunction with a unique harmonic form $F$ with $F^+ = \Omega_g$.
\end{main}

On a compact K\"ahler surface, one could therefore ask analogously if extremal K\"ahler metrics are \emph{absolute minimizers} of the Weyl functional on $\mathcal G_\Omega$ where $\Omega$ is the K\"ahler class represented by the extremal K\"ahler metric. It turns out that this is not the case:

\begin{main}\label{thmC}
For any given $\Omega\in H_{\mathrm dR}^2(M,\mathbb R)$ on K\"ahler-type smooth 4-manifolds  $\CP_2\sharp\overline{\CP}_2$ or $\CP_1\times\CP_1$ the extremal K\"ahler metrics in $\mathcal G_\Omega$ (with respect to some complex structure) are not necessarily minimizers of the Weyl functional restricted to $\mathcal G_\Omega$.
\end{main}

Theorem \ref{thmA} and \ref{thmB} are proved in Section 3 in Propositions \ref{prop1}, \ref{prop1b}, \ref{prop2}. Theorem \ref{thmC} is a consequence of the discussion in Section 4.

Recall that, given a compact complex manifold $(M,J)$ with a K\"ahler class $\Omega$ (i.e. $\Omega$ is represented by a K\"ahler form), an \emph{extremal K\"ahler metric} is, by definition, the critical point of the action
\begin{align}\label{calabi}
\Omega^+ & \longrightarrow\mathbb R\\
\nonumber \omega & \longmapsto\int_M s_{\omega}^2 d\mu_{\omega}
\end{align}
where $\Omega^+$ stands for the space of \emph{K\"ahler} forms in the de Rham class $\Omega$. This notion of extremal metrics was introduced by Calabi \cite{C} in an attempt to show existence of constant scalar curvature K\"ahler metrics on compact complex manifolds. The Euler-Lagrange equations of this action are given by $\frac{\partial^2 s}{\partial \bar{z}_i \partial \bar{z}_j} = 0$. In particular, every constant scalar curvature K\"ahler metric is extremal. The converse, however, is not true: For any given K\"ahler class on Hirzebruch srufaces $\mathbb F_k = \mathbb P(\mathcal O (-k) \oplus \mathcal O)$, Calabi constructed explicit extremal K\"ahler metrics in that class. However, the first Hirzebruch surface $\mathbb F_1 \approx \CP_2\sharp\overline{\CP}_2$ cannot admit constant scalar curvature K\"ahler metric by Matsushima--Lichnerowicz theorem because the maximal compact subgroup Lie group of automorphisms is not reductive \cite{C}.

By computing the second variation of this action at a critical metric, Calabi was able to show that extremal K\"ahler metrics are local minimizers \cite{C}. Indeed, they turn out to be \emph{global} minimizers as proven recently by Donaldson and Chen \cite{D2}, \cite{Chen}. As we will discuss in Section 2, on compact complex surfaces, LeBrun showed that the constant scalar curvature K\"ahler metrics remain to be global minimizers of the action  \eqref{calabi} if we extend the domain from $\Omega^+$ to $\mathcal G_\Omega$, provided $c_1\cdot\Omega\leq 0$. Note that $s^2/24 = |W^+|^2$ for any K\"ahler metric. However, we will show in Section 4 that the extremal K\"ahler metrics are \emph{not} necessarily global minimizers of the Weyl energy functional $g\mapsto \int |W^+|^2 d\mu$ on $\mathcal G_\Omega$.

\section{Einstein--Maxwell Equations}

This section summarizes some of the results in \cite{L08}.

Recall that the Euler--Lagrange equations of the action $g\mapsto \int_M s_g d\mu_g$, where $g$ is allowed to vary over all \textit{unit volume} metrics are precisely $\mathring{r}=0$ (i.e. Einstein metrics). Also, from Hodge theory, the Euler--Lagrange equations of the action $F\mapsto\int_M |F|^2_g d\mu_g$, where $g$ is fixed but $F$ is varying over all closed 2-forms in a fixed de Rham class $[F]\in H_{dR}^2(M,\mathbb R)$ are the Laplace equation $\Delta F=0$. Therefore, the Einstein-Maxwell equations are precisely the Euler-Lagrange equations of the joint action $(g,F)\mapsto \int_M s_g + |F|^2_g d\mu_g$ where $g$ is varying over unit volume Riemannian metrics and $F$ is varying over a fixed de Rham class.

If we restrict the first action to the conformal class of a critical metric, we get the Einstein-Hilbert action whose critical points are well known to have constant scalar curvature \cite{Yamabe1960}. Thus, any Einstein-Maxwell metric is of constant scalar curvature. Conversely, C. LeBrun observed the following remarkable fact:

\begin{prop}[LeBrun]\label{LEBprop1}
Suppose that $(M^4, g, J)$ is a K\"ahler surface with K\"ahler form $\omega=g(J\cdot,\cdot)$ and Ricci form $\rho=r(J\cdot,\cdot)$. If $g$ is constant scalar curvature K\"ahler, then $g$ satisfies Einstein-Maxwell equations with $F=\omega + \frac{1}{2} \mathring{\rho}$, where $\mathring{\rho} = \mathring{r}(J\cdot,\cdot)$ is the primitive part of the Ricci form $\rho$ of $g$.
\end{prop}

Recall that constant scalar curvature K\"ahler metrics are in particular extremal K\"ahler, which are critical points of $L^2$-norm of scalar curvature $g\mapsto\int_M s_g^2 d\mu_g$, where $g$ is varying over K\"ahler metrics on a fixed K\"ahler class $\Omega\in H_{dR}^2 (M,\mathbb R)$. C. LeBrun generalized the notion of a K\"ahler class and Calabi problem for a K\"ahler surface to the Riemannian setting, where, a priori, there may not be a complex structure at all. The generalization is as follows:

Let $M$ be a smooth 4-manifold; and let $\Omega$ be a de Rham class as above. By Hodge theory, we know that any Riemannian metric $g$ gives a unique harmonic representative $\Omega_g$ of $\Omega$. If $\Omega_g$ is self-dual, $g$ is called an $\Omega$-adapted metric. The space of all $\Omega$-adapted metrics is denoted by $\mathcal G_\Omega$; i.e. $\mathcal G_\Omega=\{g:*\Omega_g =\Omega_g\}$.

Observe that if $M$ is a complex surface and $\Omega$ is a K\"ahler class, then $\mathcal G_\Omega$ contains all K\"ahler metrics in $\Omega$, because any K\"ahler form is self-dual. In this sense, $\mathcal G_\Omega$ is a Riemannian generalization of a K\"ahler class. Also note that if $g\in \mathcal G_\Omega$, so is $\tilde g=ug\in \mathcal G_\Omega$ since Hodge $*$-operator is unchanged under conformal changes of the metric.

Now, as in the Calabi problem, C. LeBrun considers the action
\begin{equation}\tag{*}
g\mapsto \int_M s_g^2 d\mu_g
\end{equation}
on $\mathcal G_\Omega$, and sees which metrics are critical points of this action:

\begin{prop}[LeBrun]\label{LEBprop2}
Critical points of (*) are either\\
(1) scalar-flat metrics (i.e. $s\equiv 0$), or\\
(2) Einstein-Maxwell metrics $g$ with $F^+ = \Omega_g$.
\end{prop}

Thus, in particular, constant scalar curvature K\"ahler metrics are critical points of (*). Moreover, they are actually \textit{minimizers} if $c_1\cdot\Omega \leq 0$.

\begin{thm}[LeBrun]\label{LEBthm}
Let $(M^4,J)$ be a compact complex surface and $\Omega$ is a K\"ahler class with $c_1\cdot \Omega \leq 0$. Then any metric $g$ in $\mathcal G_\Omega$ satisfies  $\int_M s^2 d\mu \geq 32\pi^2 \frac{(c_1\cdot\Omega)^2}{\Omega\cdot\Omega}$, and equality holds iff $g$ is constant scalar curvature K\"ahler.
\end{thm}

Another observation of C. LeBrun is that any compact smooth 4-manifold of K\"ahler type admits a solution of \eqref{eq1}. This follows from Shu's result \cite{S} , which says that such $4$-manifolds admit a constant scalar curvature K\"ahler metrics unless they are diffeomorphic to $\mathbb CP_2\sharp\overline{\mathbb CP}_2$ or $\mathbb CP_2\sharp 2 \overline{\mathbb CP}_2$. However, both of these manifolds admit Einstein metrics (Page metric \cite{P} and Chen--LeBrun--Weber metric \cite{CLW}) which are automatically Einstein-Maxwell with $F=0$.

\section{Bach--Merkulov Equations}

In this section we will state and prove analogues of LeBrun's results stated in section 2 for Bach--Merkulov equations.

First we start by observing the following proposition which shows that Bach--Merkulov equations possess an interesting family of solutions.

\begin{prop}\label{prop1}
Let $g$ be an extremal K\"ahler metric. Then $(g,F)$ satisfies \eqref{eq2} where $F= \omega + \frac{1}{2} \psi$ where $\psi = B(J\cdot,\cdot)$. Hence any metric conformal to an extremal K\"ahler metric is a solution of \eqref{eq2}.
\end{prop}

\proof The proof is similar to the one of Proposition \ref{LEBprop1}. First, observe that $[F\circ F]_\circ = 2F^+ \circ F^-$ where $F^+$ and $F^-$ are the self-dual and anti-self-dual part of $F$, respectively. Since $g$ is K\"ahler, $\omega$ is a self-dual harmonic 2-form. Moreover, since $g$ is \textit{extremal}, $\psi = B(J\cdot,\cdot)$ is an anti-self-dual harmonic 2-form (see \cite{CLW}). Thus, setting $F^+ =\omega$ and $F^- = \frac{\psi}{2}$, we see that $2 F^+ \circ F^- = {\omega_i}^s \psi_{sj} = \psi(J\cdot,\cdot) = -B$. Thus we get $B+[F\circ F]_\circ = 0$. Moreover, $F$ is harmonic since both $F^+$ and $F^-$ are so. Therefore, $(g,F)$ is a solution of Bach--Merkulov equations. \hfill $\square$

More explicitly, if $g$ is extremal K\"ahler, then the Bach tensor can be re-written in the form
$$B =\frac{1}{12}(s\mathring r + 2 \Hess_\circ (s))$$
and therefore $\psi = \frac{1}{12}[s\rho+ i\partial\bar\partial s]_\circ$ where $[\,\cdot\,]_\circ$ stands for the primitive part of a (1,1)-form (see \cite{CLW}). In particular, if the extremal K\"ahler metric turns out to have \textit{non-zero constant scalar curvature}, then $\psi$ simplifies to $\frac{s}{12}\mathring \rho$. So we see that the solution of Proposition \eqref{prop1} becomes $(g,F=\omega+ \frac{s}{24}\mathring \rho)$ which is quite similar to LeBrun's solution to Einstein--Maxwell equations.

Proposition \ref{prop1} together with Shu's result implies the following:

\begin{prop}\label{prop1b}
Let $M$ be the underlying $4$-manifold of any compact complex surface of K\"ahler type. Then $M$ admits a solution $(g,F)$ of Bach--Merkulov equations.
\end{prop}

Next, we will prove the analogue of Proposition \ref{LEBprop2} for Bach--Merkulov equations:

\begin{prop}\label{prop2}
An $\Omega$-adapted metric $g$ is a critical point of the restriction of Weyl functional to $\mathcal G_\Omega$ iff $g$ is a solution of Bach--Merkulov equations in conjunction with a unique harmonic form $F$ with $F^+ = \Omega_g$.
\end{prop}
\proof The proof is similar to the one of Proposition \ref{LEBprop2}. Let $g_t = g+th+O(t^2)$ be a variation of a metric $g$ in $\mathcal G_\Omega$. Donaldson showed that the tangent space $T\mathcal G_\Omega$ is precisely the $L^2$-orthogonal complement of $\{\Omega_g\circ \varphi : \varphi\in \mathcal H_g^- \}$ in $\Gamma(\bigodot^2 T^*M)$. Thus, in our case, $h$ can be taken such that $\int_M \langle h, \Omega_g \circ \varphi\rangle d\mu_g = 0$ for all $\varphi \in \mathcal H_g$.

The first variation of the Weyl functional is given by (\cite{B}, \cite{Bes})
$$
\frac{d}{dt}\int_M \|W\|^2 d\mu_{g_t} = \int_M h^{ij}B_{ij} d\mu_g = \int_M \langle h,B\rangle d\mu_g.
$$
Thus, $g$ is a critical point iff $h$ is $L^2$-orthogonal to $B$. By Donaldson's result, this implies that $B=\Omega_g\circ\varphi$ for some $\varphi \in \mathcal H_g^-$. So, taking $F^+ =\Omega_g$ and $F^- = -\varphi$, we see that $g$ satisfies \eqref{eq2}.

Conversely, if $g$ is an $\Omega$-adapted solution of \eqref{eq2} with $F^+=\Omega_g$, then $\int \langle h,B\rangle d\mu = \int \langle h, F^+ \circ F^-\rangle d\mu = 0$ for any variation $h$ as above. Thus, by Donaldson, $g$ is a critical point. \hfill $\square$

In particular, extremal K\"ahler metrics are also critical points of this functional. The natural question to ask is whether they are global minimizers in $\mathcal G_\Omega$. In the next section, we will show that the answer to this question is negative: the analogue of Theorem \ref{LEBthm} does not hold for Bach--Merkulov equations.

\section{Example: Hirzebruch Surfaces}

In this section we will show that extremal K\"ahler metrics adapted to a fixed cohomology class $\Omega$ do no necessarily have the same Weyl energy. We will illustrate this fact on Hirzebruch surfaces, by showing existence of two closed forms in $\Omega$ which are extremal K\"ahler with respect to different complex structures. Using the formula in \cite{HS97} it will turn out that the Weyl energy of the corresponding extremal K\"ahler metrics are different.

Recall that the $k$-th Hirzebruch surface $\mathbb F_k$ is defined as the projectivization of the rank-2 complex vector bundle $\mathcal O (-k) \oplus \mathcal O$ over $\mathbb C P_1$ (see \cite{Bea} and \cite{BPV} for details). $\mathbb F_k$ is diffeomorphic to $S^2\times S^2$ if $k$ is even, and to $\mathbb CP_2\sharp \overline{\mathbb CP_2}$ if $k$ is odd $\cite{Hirz}$. They are, however, all biholomorphically distinct as complex surfaces (see \cite{Hirz}). They are simply connected; they have second Betti number $b_2 (\mathbb F_k) =2$ and Euler characteristic $\chi (\mathbb F_k)$.

For the generators of the homology of $\mathbb F_k$ we will take the fiber $F$ and the image of the section $\{z\mapsto [0:z]\} : \mathbb C P_1 \rightarrow \mathbb P (\mathcal O(-k) \oplus \mathcal O) = \mathbb F_k$, which we will denote by $C_k$. Note that $F\cdot F =0$, $F\cdot C_k = 1$ and $C_k\cdot C_k = -k$ so that in this basis the intersection pairing becomes $$ \left(\begin{array}{rr} 0 &1 \\ 1 &-k\end{array}\right).$$

Let the Poincar\'e dual of $C_k$ and $F$ be $\mathfrak c_k$ and $\mathfrak f$ respectively. Then any de Rham class $\Omega\in H^2_{dR}(\mathbb F_k, \mathbb R)$ can be written as $\Omega = p\mathfrak c_k + q \mathfrak f$ for some $p,q\in\mathbb R$. If $k$ and $n$ are two positive integers of same parity, then $\mathbb F_k$ and $\mathbb F_n$ are diffeomorphic; so we can represent $\Omega$ with respect to the basis $\{\mathfrak c_n, \mathfrak f\}$. The following lemma gives the change of basis formula:

\begin{lem}\label{lemtransf}
We have
$$\mathfrak c_k = \mathfrak c_n + \frac{n-k}{2} \mathfrak f.$$
Therefore,
$$\Omega = p\mathfrak c_k + q\mathfrak f = p \mathfrak c_n + \tilde q\mathfrak f$$
where $\tilde q = p \frac{n-k}{2} + q$
\end{lem}
\proof Let
\begin{equation}\label{eqckcnf}
C_k = s C_n + t F
\end{equation}
 for some constants $s,t$. Take the intersection of both sides with $F$:
$$ C_k\cdot F = s C_n \cdot F + t F\cdot F. $$
Since $C_k\cdot F = C_n\cdot F =1$ and $F\cdot F = 0$, we have $s=0$. On the other hand, take the self intersection of both sides:
\begin{eqnarray}\nonumber
C_k \cdot C_k &=& (C_n + t F)\cdot (C_n + t F)\\ \nonumber
-k &=& -n + 2t.
\end{eqnarray}
Therefore $t= \frac{n-k}{2}$. Taking the Poincar\'e dual of \eqref{eqckcnf} proves the first equality. The second equality follows immediately from this. \hfill$\square$

Note also that $\Omega = p\mathfrak  c_k + q \mathfrak f \in H^2_{dR} (\mathbb F_k, \mathbb R)$ is a K\"ahler form iff $\Omega\cdot C_k>0$ and $\Omega \cdot F>0$, that is, iff $p>0$ and $q>kp$. Now we can deduce when the same de Rham class $\Omega$ is a K\"ahler class in $\mathbb F_n$, where $n$ and $k$ have the same parity. Let $J_k$ denote the complex structure of the complex surface $\mathbb F_k$.

\begin{lem}
A K\"ahler class $\Omega = p \mathfrak c_k +q \mathfrak f$ in $\mathbb F_k$ is a K\"ahler class in $\mathbb F_n$ iff $n< 2\frac{q}{p}-k$. In particular, $\Omega$ is K\"ahler with respect to only \emph{finitely many} $J_n$'s.
\end{lem}

\proof $\Omega = p \mathfrak c_k +q \mathfrak f$ is K\"ahler with respect to $J_k$ iff $p>0$ and $q>kp$. By Lemma \ref{lemtransf}, $\Omega = p\mathfrak c_n + \left(p\frac{n-k}{2} + q\right)$. Now, by the previous paragraph, this class is K\"ahler with respect to $J_n$ iff $p>0$ and $p\frac{n-k}{2} + q>np$. The second inequality is the same as $n<2\frac{q}{p} - k$. Notice that $2\frac{q}{p} - k$ is positive since $\Omega = p \mathfrak c_k +q \mathfrak f$ is assumed to be K\"ahler since $q>pk$. Hence there are only finitely many possibilities for $n$ so that $\Omega$ remains K\"ahler with respect to $J_n$. \hfill$\square$

So there are de Rham classes on smooth 4-manifolds $S^2\times S^2$ or $\mathbb CP_2 \sharp \overline{\mathbb C P_2}$ which are K\"ahler with respect to different complex structures. However, Calabi \cite{C} showed that every K\"ahler class on a Hirzebruch surface is represented by an \emph{extremal} K\"ahler metric.

So, with our previous notation all Riemanninian metrics $g$ whose K\"ahler form $\omega=g(J_k\cdot,\cdot)$ with respect to any of the complex structures $J_k$ are in $\mathcal G_\Omega$. Thus, we have essentially distinct extremal K\"ahler metrics in $\mathcal G_\Omega$. Each of those metrics are critical points of the restriction of the Weyl functional to $\mathcal G_\Omega$.

Next, we will show that this the Weyl energy levels of those metrics are different. First note that for K\"ahler metrics we have $|W|^2 = \frac{s^2}{24}$. Therefore the Weyl energy of a K\"ahler metric is equal to its Calabi energy up to a constant of 24. Hwang\&Simanca \cite{HS97} gave the following formula for the Calabi energy of an extremal metric in a K\"ahler class on the Hirzebruch surface $\mathbb F_k$.

\begin{prop}[Hwang\&Simanca]
The Calabi energy of the extremal K\"ahler metric in the class $\Omega= 4\pi \mathfrak c_k + 2\pi (a+k) \mathfrak f$ is given as:
\begin{equation}\label{eqhwangsimanca}
\tilde {\mathcal C} (a,k) := 12\pi \frac{a^3 +4a^2+(4+k^2)a -4k^2}{3a^2 -k^2}.
\end{equation}
\end{prop}
Note that the Calabi energy and the Weyl energy are scale-invariant in dimension four. Therefore by appropriate scaling we see that the Calabi energy of the extremal K\"ahler metric in $\Omega=pc_k + qf$ is given by
\begin{equation}
\mathcal C(p,q,k) := \tilde {\mathcal C}(2\frac{q}{p}-k,k) = 12\pi \frac{(2\frac{q}{p}-k)^3 +4(2\frac{q}{p}-k)^2+(4+k^2)(2\frac{q}{p}-k) -4k^2}{3(2\frac{q}{p}-k)^2 -k^2}.
\end{equation}

We therefore see that the extremal K\"ahler metrics with respect to different complex structures in $\Omega$ have different energy, i.e. $\mathcal C(p,q,k) \neq \mathcal C(p,p \frac{n-k}{2}+q, n)$ in general.

This shows that the analogue of of Theorem \ref{LEBthm} cannot hold for the Bach--Merkulov equations.




\noindent
{\bf Acknowledgement. } I would like to thank my advisor Claude LeBrun for suggesting the problem and for his help, guidance and encouragement.

\vspace{1in}

\noindent
{\sc Author's address:}

\medskip

 \noindent
{Mathematics Department, SUNY, Stony Brook, NY 11794, USA
}

\bigskip

\noindent
{\sc Author's e-mail:}

\medskip

 \noindent
{\tt caner@math.sunysb.edu
}

\end{document}